\documentclass[preprint,authoryear]{elsarticle}

\biboptions{sort&compress}

\usepackage{graphicx}      
\usepackage{natbib}        
\usepackage{amssymb}
\usepackage{amsmath}
\usepackage{enumerate}
\usepackage{placeins}
\usepackage{hyphenat}
\usepackage{epstopdf}

\newtheorem{thm}{Theorem}

\newtheorem{corol}{Corollary}
\newtheorem{assum}{Assumption}
\newproof{pf}{Proof}

\begin{document}
\bibliographystyle{model5-names}
\pagestyle{plain}
\setcounter{page}{0}
\pagenumbering{arabic}
\begin{frontmatter}

\title{On the Equivalence Between the Modifier-Adaptation and Trust-Region Frameworks} 

\author{Gene A. Bunin}
\ead{gene.a.bunin@ccapprox.info}

\begin{abstract}                

In this short note, the recently popular modifier-adaptation framework for real-time optimization is discussed in tandem with the well-developed trust-region framework of numerical optimization, and it is shown that the basic version of the former is a simplification of the latter when the problem is unconstrained. This relation is then exploited to propose a globally convergent modifier-adaptation algorithm using already developed trust-region theory. Cases when the two may not be equivalent and extensions to constrained problems are also discussed.

\;


\end{abstract}

\end{frontmatter}

\section{The Real-Time Optimization Problem}
\label{sec:intro}

In the process systems engineering community, the basic idea of most real-time optimization (RTO) schemes consists in finding a set of optimal operating conditions -- often steady-state setpoints in a multilayer hierarchical scheme -- that minimize (resp., maximize) the steady-state cost (resp., profit) of some given plant subject to constraints \citep{Brdys2005}. While models of the process being optimized are often available, it is generally the case that they are either inaccurate and/or incomplete, which motivates the data-driven ``real-time'' element of RTO, thereby forcing the optimization algorithm to use the measurements obtained from the process as feedback to modify the provided setpoints so as to ultimately reject the model uncertainty and converge to the optimal conditions of the plant. 

A fairly general mathematical formulation of this problem that suffices for many practical cases is as follows:

\begin{equation}\label{eq:mainprobcons}
\begin{array}{rl}
\mathop {{\rm{minimize}}}\limits_{\bf{u}} & \phi_{p} ({\bf{u}}) \\
{\rm subject\;to} & g_{p,j}({\bf u}) \leq 0,  \;\; j = 1,...,n_g,
\end{array}
\end{equation}

\noindent where ${\bf u} \in \mathbb{R}^{n_u}$ denote the decision variables, or the ``inputs'', of the problem, while the functions $\phi, g : \mathbb{R}^{n_u} \rightarrow \mathbb{R}$ denote the cost and constraints, respectively. The subscript $p$ (for ``plant'') is used to indicate that the function corresponds to an experimental relationship that is not perfectly known and may only be approximated by a model, which we will mark with the subscript $\hat p$ (e.g., $\phi_{\hat p}$ being the model approximation of $\phi_{p}$). In the simplest terms, the goal of an RTO algorithm is to solve Problem (\ref{eq:mainprobcons}) by iterative experimentation, generating a sequence of steady-state ${\bf u}$ values that converges to the plant optimum.

For the majority of this document, we will not focus on Problem (\ref{eq:mainprobcons}) but on the unconstrained case

\begin{equation}\label{eq:mainprob}
\mathop {{\rm{minimize}}}\limits_{\bf{u}} \;\;  \phi_p ({\bf{u}})\;\;,
\end{equation}

\noindent as this is sufficient to convey the main message. We will, however, return to Problem (\ref{eq:mainprobcons}) in the end in passing, providing references to works where it is discussed properly and in much greater detail. 

\section{Review of the Modifier-Adaptation Framework}
\label{sec:review}

An approach to solving (\ref{eq:mainprobcons}) that has recently gained popularity in the research community is that of \emph{modifier adaptation}, which originally dates back to the work of \cite{Roberts1978} and owes its numerous refinements and fundamental ideas to the \emph{ISOPE} (``iterative setpoint optimization and parameter estimation'') framework \citep{Brdys2005}. Recent works by \cite{Gao2005}, \cite{Chachuat2009}, and \cite{Marchetti2009a} have given the approach its modern form by accounting for plant-model mismatch in both the cost and constraints. A number of works in the past few years have also considered various particular aspects of the framework, such as mathematical reformulations to ease or better accommodate particular problem types \citep{Francois2013,Serralunga2013,Costello2013}, important implementation aspects \citep{Marchetti2010,Rodger2010,Bunin2012}, and  major theoretical issues like feasibility \citep{Bunin2011,Navia2012} and global convergence \citep{Faulwasser2014}.

The basic philosophy of modifier adaptation lies in applying local corrections to an inherently incorrect model at each RTO iteration $k$, and solving this corrected version to obtain the following iterate at $k+1$. For the unconstrained case, this would lead to the following update:

\begin{equation}\label{eq:newu}
\begin{array}{rl}
{\bf u}_{k+1} \in {\rm arg}\; \mathop {{\rm{minimize}}}\limits_{\bf{u}} & \phi_{\hat p} ({\bf{u}}) + {\boldsymbol \lambda}_k^T {\bf u}\;\;,
\end{array}
\end{equation}

\noindent with the \emph{modifiers} ${\boldsymbol \lambda}_k$, defined as

\begin{equation}\label{eq:lambda}
{\boldsymbol \lambda}_k := \nabla \phi_p ({\bf u}_k) - \nabla \phi_{\hat p} ({\bf u}_k),
\end{equation}

\noindent serving to ensure that the plant and corrected model have matching first derivatives at the current iterate ${\bf u}_k$.

Placing this into algorithmic form yields the following basic implementation.

\newpage
\noindent {\bf Algorithm 1 (Basic Modifier-Adaptation Algorithm)}
\;

\begin{enumerate}
\item (Initialization) The initial point, ${\bf u}_0$, is provided. Set $k := 0$.
\item (Modifier Computation) Compute the modifiers ${\boldsymbol \lambda}_k$ according to (\ref{eq:lambda}).
\item (New Input Calculation) Obtain ${\bf u}_{k+1}$ by solving Problem (\ref{eq:newu}) and apply this set of inputs to the plant.
\item (Iterate) Set $k := k+1$ and return to Step 2.  
\end{enumerate}

The key oft-stated motivation for applying this algorithm is the following upon-convergence guarantee.

\begin{thm}[First-Order Critical Point Upon Convergence]
Assume that the minimization of (\ref{eq:newu}) always yields a first-order critical point of the modified objective function $\phi_{\hat p} ({\bf{u}}) + {\boldsymbol \lambda}_k^T {\bf u}$ and that Algorithm 1 has converged to a fixed point ${\bf u}_\infty$. It follows that ${\bf u}_\infty$ is a first-order critical point of $\phi_{p}$.
\end{thm}
\begin{pf}
The result follows immediately from the fact that a first-order critical point for an unconstrained problem is defined entirely by the function's derivatives at that point. As these must match for the modified model and the plant at any iterate, including ${\bf u}_\infty$, it follows that finding a first-order critical point for the modified function implies finding one for the plant. \qed
\end{pf}

\section{The Basic Trust-Region Algorithm}
\label{sec:TRalgo}

A theoretically rigorous approach for iteratively minimizing a nonlinear function in the mathematical optimization context is that of trust-region methods. In this section, we will consider what attempting to solve Problem (\ref{eq:mainprob}) in this framework would entail.

Let us start by stating the basic trust-region algorithm for solving (\ref{eq:mainprob}). This is essentially the algorithm provided in the well-known monograph on trust-region methods \cite[Ch. 6]{Conn2000} but with a few additional simplifications and some notational changes. Namely, we use the 2-norm instead of the general $p$-norm and explicitly distinguish between the reference iterates, ${\bf u}_k^*$, and the iterates applied to the plant, ${\bf u}_k$.

\;
\;
\noindent {\bf Algorithm 2 (Basic Trust-Region Algorithm)}
\;

\begin{enumerate}
\item (Initialization) The initial point, ${\bf u}_0$, and initial trust-region radius, $\Delta_0 > 0$, are provided, together with the constants $\eta_1$, $\eta_2$, $\gamma_1$, and $\gamma_2$ satisfying $0 < \eta_1 \leq \eta_2 < 1$ and $0 < \gamma_1 \leq \gamma_2 < 1$. Set $k := 0$, ${\bf u}_0^* := {\bf u}_0$, and apply ${\bf u}_0$ to the plant to obtain $\phi_p ({\bf u}_0^*) = \phi_p ({\bf u}_0)$.

\item (Model Construction) Construct the model $m_k$, which is an approximation of $\phi_p$ over the trust region $\mathcal{B}({\bf u}_k^*,\Delta_k)$, i.e., over a Euclidean ball of radius $\Delta_k$ centered at ${\bf u}_k^*$.
\item (New Input Candidate Calculation) Compute a candidate point ${\bf u}_{k+1} \in \mathcal{B}({\bf u}_k^*,\Delta_k)$ that ``sufficiently reduces the model'' $m_k$.
\item (Acceptance of the Candidate Point) Apply ${\bf u}_{k+1}$ to the plant and evaluate $\phi_p ({\bf u}_{k+1})$. Define:

\begin{equation}\label{eq:rho}
\rho_k := \frac{\phi_p({\bf u}_k^*) - \phi_p({\bf u}_{k+1})}{m_k({\bf u}_k^*) - m_k({\bf u}_{k+1})}.
\end{equation}
 
\noindent If $\rho_k \geq \eta_1$, then set ${\bf u}_{k+1}^* := {\bf u}_{k+1}$. Otherwise, set ${\bf u}_{k+1}^* := {\bf u}_k^*$.

\item (Trust-Region Radius Update) Set $\Delta_{k+1}$ such that

\begin{equation}\label{eq:trupdate}
\Delta_{k+1} \in \left\{ 
\begin{array}{ll}
[ \Delta_k, \infty) & {\rm if}\; \rho_k \geq \eta_2,  \\
\left[ \gamma_2 \Delta_k, \Delta_k \right]  & {\rm if}\; \rho_k \in [\eta_1, \eta_2 ), \\
\left[ \gamma_1 \Delta_k, \gamma_2 \Delta_k \right] & {\rm if}\; \rho_k < \eta_1.
\end{array} \right .
\end{equation}

\item (Iterate) Set $k := k+1$ and return to Step 2.

\end{enumerate}

Let us now state the assumptions sufficient to prove the global convergence of Algorithm 2 to a first-order critical point \citep{Conn2000}. The following are assumed about the nature of the plant:

\begin{assum}
$\phi_p$ is $\mathcal{C}^2$ (twice continuously differentiable) on $\mathbb{R}^{n_u}$.
\end{assum}

\begin{assum}
$\phi_p$ is lower-bounded on $\mathbb{R}^{n_u}$.
\end{assum}

\begin{assum}
The Hessian of $\phi_p$ is upper-bounded on $\mathbb{R}^{n_u}$.
\end{assum}

\noindent As mentioned in \cite{Conn2000}, Assumption 3 is often too strong and could actually be restricted to the subspace of $\mathbb{R}^{n_u}$ where the iterates lie. However, as this subspace is not known \emph{a priori}, $\mathbb{R}^{n_u}$ is used for notational convenience.

The following assumptions are made on the model:

\begin{assum}
For all $k$, $m_k$ is $\mathcal{C}^2$ over $\mathcal{B}({\bf u}_k^*, \Delta_k)$.
\end{assum}

\begin{assum}
$m_k$ matches $\phi_p$ locally to first order at every $k$, i.e.:

\begin{equation}\label{eq:match1}
m_k({\bf u}_k^*) = \phi_p ({\bf u}_k^*),
\end{equation}

\begin{equation}\label{eq:match2}
\nabla m_k({\bf u}_k^*) = \nabla \phi_p ({\bf u}_k^*).
\end{equation}

\end{assum}

Finally, one requires the following assumption on the algorithm used to solve the trust-region subproblem with regard to its ability to achieve ``sufficient reduction'' in the model:

\begin{assum}
There exists a constant $\kappa \in (0,1)$ such that for all $k$:

\begin{equation}\label{eq:decrease}
m_k ({\bf u}_k^*) - m_k ({\bf u}_{k+1}) \geq \kappa \| \nabla m_k ({\bf u}_k^*) \| \mathop {\min} \left[ \frac{\| \nabla m_k ({\bf u}_k^*) \|}{\beta_k}, \Delta_k \right],
\end{equation}

\noindent with $\beta_k > 1$ a finite constant.

\end{assum}

One may then state the following.

\begin{thm}[Global Convergence to a First-Order Critical Point]
If Assumptions 1-6 are satisfied, it then follows that the iterates generated by Algorithm 2 converge asymptotically to a first-order critical point, i.e.:

\begin{equation}\label{eq:globconv}
\mathop {\lim} \limits_{k \rightarrow \infty} \| \nabla \phi_p ({\bf u}_k^*) \| = 0.
\end{equation}

\end{thm}
\begin{pf}
The reader is referred to Theorem 6.4.6 in \cite{Conn2000}. Note that we have, for simplicity, used a slightly stronger assumption and have assumed that $m_k$ is $\mathcal{C}^2$ over $\mathcal{B}({\bf u}_k^*, \Delta_k)$. The two assumptions made by \cite{Conn2000} -- namely, that over $\mathcal{B}({\bf u}_k^*, \Delta_k)$ the model  $m_k$ is twice differentiable and that its Hessian is bounded -- are implied by the single $\mathcal{C}^2$ assumption here. \qed
\end{pf}

\section{Equivalence and a Globally Convergent Modifier-Adaptation Scheme}
\label{sec:equi}

Both the modifier-adaptation and trust-region algorithms seek to minimize $\phi_p$ by iteratively optimizing a local approximation of $\phi_p$ around each ${\bf u}_{k}^*$. The key differences between the two may be summarized as follows:

\begin{enumerate}
\item The model $m_k({\bf u}) = \phi_{\hat p} ({\bf u}) + {\boldsymbol \lambda}_k^T {\bf u}$ used by modifier-adaptation enforces, by construction, (\ref{eq:match2}) but not (\ref{eq:match1}). The standard trust-region algorithm usually enforces both as this is required by the convergence proof.
\item The modifier-adaptation subproblem (\ref{eq:newu}) considers the whole input space while the trust-region subproblem limits its search to the ball $\mathcal{B}({\bf u}_k^*, \Delta_k)$.
\item The concept of a ``reference point'' is absent in the basic modifier-adaptation algorithm, as the computed ${\bf u}_{k+1}$ is always used as the reference with respect to which the model is corrected at the subsequent iteration. In the trust-region scheme, the model is always built with respect to the latest ``successful'' iterate for which a sufficient decrease in the plant cost function value has been achieved.
\end{enumerate}

The first difference is actually of no practical consequence -- as discussed later (see Corollary 1), one could always use a model that satisfies both (\ref{eq:match1}) and (\ref{eq:match2}) without changing the iterates generated by the modifier-adaptation algorithm. The second and third differences, however, are important and may aid in explaining why no globally convergent version of Algorithm 1 has been derived to date. Without the use of a reference point, it is difficult to ensure the stability of the algorithm, since any progress made may always be undone by a single bad iteration. Optimizing with respect to the best known point effectively prevents bad iterations from having any lasting effect on convergence, but is not sufficient to guarantee the existence of a good iteration. For this, one needs the guarantee that the model used by the algorithm become sufficiently good under certain conditions. Since the model is only good locally and to first-order, the natural approach, and the one pursued in trust-region methods, is to shrink the search space until this approximation is good enough to generate a successful iterate. By considering the entire input space, the modifier-adaptation algorithm may generate iterates in portions of the input space that are not accurately modeled, and so it should not be  surprising that the guarantee of successful iterates is absent in this algorithm.

Note that all of these differences are of the same nature, in that they are all things that are present in the trust-region framework but absent in modifier adaptation. In fact, if we were to enforce that ${\bf u}_{k+1}^* := {\bf u}_{k+1}$ always and let $\Delta_0 \rightarrow \infty$ and $\gamma_1 \rightarrow 1$ (i.e., remove the trust-region restriction) in Algorithm 2, we would essentially end up with Algorithm 1. Considering things from this perspective, let us now avoid these simplifications and propose the following modifier-adaptation scheme.

\;
\;
\noindent {\bf Algorithm 3 (Trust-Region Supplemented Modifier-Adaptation Algorithm)}
\;

\begin{enumerate}
\item (Initialization) Identical to Step 1 of Algorithm 2.
\item (Modifier Computation) Compute the modifiers ${\boldsymbol \lambda}_k := \nabla \phi_p ({\bf u}_k^*) - \nabla \phi_{\hat p} ({\bf u}_k^*)$.
\item (Model Correction) Construct the model $m_k({\bf u}) := \phi_{\hat p} ({\bf u}) + {\boldsymbol \lambda}_k^T {\bf u}$.
\item (New Input Candidate Calculation) Compute a candidate point ${\bf u}_{k+1}$ by approximately solving the problem

\begin{equation}\label{eq:newu2}
\begin{array}{rl}
{\bf u}_{k+1} \in {\rm arg}\; \mathop {{\rm{minimize}}}\limits_{{\bf u} \in \mathcal{B}({\bf u}_k^*, \Delta_k)} & m_{k} ({\bf{u}})\;\;.
\end{array}
\end{equation}

\noindent Furthermore, compute the Cauchy point, ${\bf u}_{k+1}^{\rm CP}$, via the line search

\begin{equation}\label{eq:newuCP}
\begin{array}{rl}
\left [{\bf u}_{k+1}^{\rm CP}, t^{\rm CP} \right] \in {\rm arg}\; \mathop {{\rm{minimize}}}\limits_{\footnotesize{\begin{array}{c} {\bf u} \in \mathcal{B}({\bf u}_k^*, \Delta_k) \\ t \geq 0 \end{array}}} & m_{k} ({\bf{u}}) \\
{\rm subject\;to} \hspace{4mm} & {\bf u} = {\bf u}_k^* - t \nabla m_k ({\bf u}_k^*).
\end{array}
\end{equation}

\noindent If $m_k ({\bf u}_{k+1}) > m_k ({\bf u}_{k+1}^{\rm CP})$, set ${\bf u}_{k+1} := {\bf u}_{k+1}^{\rm CP}$.

\item (Acceptance of the Candidate Point) Identical to Step 4 of Algorithm 2.
\item (Trust-Region Radius Update) Identical to Step 5 of Algorithm 2.
\item (Iterate) Set $k := k+1$ and return to Step 2.

\end{enumerate}

Prior to proving the global convergence of Algorithm 3, we modify Assumption 4 to make it more direct.

\begin{assum}
For all $k$, $\phi_{\hat p}$ is $\mathcal{C}^2$ over $\mathcal{B}({\bf u}_k^*, \Delta_k)$.
\end{assum}

The following key result follows.

\begin{corol}{\bf (Global Convergence to a First-Order Critical Point for Modifier Adaptation)}
If Assumptions 1-3 and 7 are satisfied, it then follows that the iterates generated by Algorithm 3 converge asymptotically to a first-order critical point, i.e.:

\begin{equation}\label{eq:globconv2}
\mathop {\lim} \limits_{k \rightarrow \infty} \| \nabla \phi_p ({\bf u}_k^*) \| = 0.
\end{equation}

\end{corol}
\begin{pf}
Algorithm 3 is special case of Algorithm 2, and so we just need to show that all of the assumptions needed for Theorem 2 are satisfied either implicitly or explicitly. As Assumptions 1-3 are made explicitly throughout, we focus on Assumptions 4-6. Since adding a linear correction term to a $\mathcal{C}^2$ function will not jeopardize the $\mathcal{C}^2$ property, making Assumption 7 implies that Assumption 4 holds. While Condition (\ref{eq:match2}) of Assumption 5 is satisfied by construction, Condition (\ref{eq:match1}) is not. However, note that we may just as easily use the model $m_k({\bf u}) := \phi_{\hat p} ({\bf u}) + [\phi_{ p} ({\bf u}_k^*) - \phi_{\hat p} ({\bf u}_k^*) ] +  {\boldsymbol \lambda}_k^T ( {\bf u} - {\bf u}_k^*)$, which satisfies both (\ref{eq:match1}) and (\ref{eq:match2}) by construction but does not influence the sequence of iterates produced by Algorithm 3 since the addition of the constant term $\phi_{ p} ({\bf u}_k^*) - \phi_{\hat p} ({\bf u}_k^*)  -  {\boldsymbol \lambda}_k^T  {\bf u}_k^*$ does not influence the computation of $\rho_k$ or ${\bf u}_{k+1}$ in any way. By sleight of hand, we may thus ``pretend'' to use the latter model and consider Assumption 5 satisfied, as the two models are equivalent with respect to the sequence of iterates generated. Finally, overriding the standard computation of ${\bf u}_{k+1}$ in Step 4 with the Cauchy point when needed ensures that Assumption 6 is met \citep[\textsection 6.3]{Conn2000}. \qed
\end{pf}

\section{Nonequivalent Cases and Practical Considerations}

As the guarantee of global convergence is a very desirable property, and as the additions to ensure it for the basic modifier-adaptation scheme are simple and algorithmic in nature, it is tempting to ask if not every modifier-adaptation scheme could be cast in a globally convergent trust-region formulation. While further research is required to give a definitive answer, a preliminary inspection seems to suggest the answer to be positive.

Perhaps of greatest interest is the question of how the prior discussion generalizes to the constrained problem (\ref{eq:mainprobcons}), since almost all problems in practice are constrained. The standard approach in trust-region methods is to cast such problems as unconstrained problems with a penalty for constraint violations included in the augmented cost function \citep[Ch. 14]{Conn2000}, and the recent work by \cite{Biegler2014}, without stating so explicitly, essentially shows how the constrained modifier-adaptation problem may be solved in the trust-region framework by exploiting this approach. While one could propose different implementation routes with regard to particular algorithmic aspects, there appears to be no reason as to why the generalization to (\ref{eq:mainprobcons}) would not come easily.

Another popular technique in modifier-adaptation schemes is to filter the modifiers \citep{Marchetti2009a,Chachuat2009,Serralunga2013} so as to not ``overcorrect'' the model, and to define them as

\begin{equation}\label{eq:lambdaf}
{\boldsymbol \lambda}_k := \alpha \left[ \nabla \phi_p ({\bf u}_k^*) - \nabla \phi_{\hat p} ({\bf u}_k^*) \right] + (1-\alpha) {\boldsymbol \lambda}_{k-1},
\end{equation}

\noindent starting from some initial values ${\boldsymbol \lambda}_{-1}$, with $\alpha \in (0,1]$ a filter gain. For $\alpha < 1$, the crucial Condition (\ref{eq:match2}) of Assumption 5 is generally not satisfied, and one thus cannot apply the same global convergence analysis to such algorithms. However, this technique of ``model filtering'' is very similar in essence to the ``models with memory'' discussed in the trust-region literature \citep[\textsection 9.5]{Conn2000}, and so it would not be surprising if the analysis of the latter were directly applicable to modifier-adaptation schemes that employed a filter.

Finally, it is important to emphasize that much of the discussion so far has focused on very idealized cases, without considering how the algorithms would behave in real application, where neither accurate function or derivative values would be available and where numerous other implementation issues could enter to complicate analysis \citep{Quelhas:12,Bunin:SCFOImp}. While recent research on trust-region methods has looked into cases with corrupted function values and derivatives \citep{Larson2012}, it is probably too early for such methods to be directly applicable to many practical real-time optimization problems. Nevertheless, there is no reason to suspect why theory developed for such problems not be equally applicable to both frameworks.

\bibliography{MATR}             

\end{document}